\documentclass[10pt]{amsart}

\usepackage[english]{babel}
\usepackage{pstricks}
\usepackage[applemac]{inputenc}
\usepackage{ae, aeguill}
\usepackage{amsmath,amsthm,amssymb,latexsym,amscd}
\usepackage[mathscr]{eucal}
\usepackage{pstricks}
\usepackage{pst-fill}
\usepackage{xy}
\xyoption{all}
\xyoption{poly}
\xyoption{ps}
\UsePSspecials{dvips}

\newif\ifpdf
\ifx\pdfoutput\undefined
\pdffalse 
\else
\pdfoutput=1 
\pdftrue
\fi

\ifpdf
\usepackage[pdftex]{graphicx}
\else
\usepackage{graphicx}
\fi

\newcommand{\ie}{\emph{i.e.}\:}

\newcommand{\SU}{\mathrm{SU}(2)}
\newcommand{\G}{\mathfrak{G}}
\newcommand{\SL}{\mathrm{SL}_2(\CC)}
\newcommand{\PSL}{\mathrm{PSL}_2(\CC)}
\newcommand{\SO}{\mathrm{SO}(3)}
\newcommand{\su}{\mathfrak{su}(2)}
\newcommand{\g}{\mathfrak{g}}
\newcommand{\sll}{\mathfrak{sl}_2(\CC)}

\newcommand{\ii}{\mathbf{i}}
\newcommand{\jj}{\mathbf{j}}

\newcommand{\kk}{\mathbf{k}}
\newcommand{\I}{\mathbf{1}}
\newcommand{\ZZ}{\mathbb{Z}}
\newcommand{\CC}{\mathbb{C}}
\newcommand{\IR}{\mathbb{R}}

\newcommand{\tangent}[2]{T_{#1} #2}

\newcommand{\bord}{\partial}
\newcommand{\lk}{\ell\mathit{k}}
\newcommand{\im}{\mathrm{im}}

\newcommand{\rk}{\mathop{\mathrm{rk}}\nolimits}
\begin{document}
\theoremstyle{plain}
\newtheorem*{MultLemma}{Multiplicativity Lemma}
\newtheorem{theorem}{Theorem}
\newtheorem*{theorem*}{Theorem}
\newtheorem*{maintheorem*}{Main Theorem}
\newtheorem{prop}[theorem]{Proposition}
\newtheorem{corollary}[theorem]{Corollary}
\newtheorem{lemma}[theorem]{Lemma}
\newtheorem{claim}[theorem]{Claim}
\newtheorem{fact}[theorem]{Fact}
\theoremstyle{definition}
\newtheorem{definition}{Definition}
\newtheorem{example}{Example}
\newtheorem*{example*}{Example}
\newtheorem{notation}{Notation}
\newtheorem*{notation*}{Notation}
\theoremstyle{remark}
\newtheorem{remark}{Remark}
\newtheorem{question}{Question}

\title{Non abelian twisted Reidemeister torsion for fibered knots}
\author{J\'er\^ome Dubois}
\address{Section de Math\'ematiques \\ Universit\'e de Gen\`eve CP 64, 2--4 Rue du Lièvre \\ CH 1211 Genève 4 (Suisse)}
\email{Jerome.Dubois@math.unige.ch}
\date{\today}
\ifpdf
\DeclareGraphicsExtensions{.pdf, .jpg, .tif}
\else
\DeclareGraphicsExtensions{.eps, .jpg}
\fi
\begin{abstract} 
	In this article, we give an explicit formula to compute the non abelian twisted sign-determined Reidemeister torsion of the exterior of a fibered knot in terms of its monodromy. As an application, we give explicit formulae for the non abelian Reidemeister torsion of torus knots and of the figure eight knot.
\end{abstract}
\maketitle
\markboth{Jérôme Dubois }{Non abelian twisted Reidemeister torsion for fibered knots}
\pagestyle{myheadings}

\section{Introduction}
\label{}

	A knot in the $3$-sphere is called \emph{fibered} if its exterior has the structure of a surface bundle over the circle. 
	For example, each torus knot is fibered and the figure eight knot is also fibered. The aim of this paper is to compute the sign-determined non abelian Reidemeister torsion defined by the author in~\cite{JDTHESE} (see also~\cite{CRAS} and~\cite{torsionvol}) for a fibered knot in $S^3$ in terms of the eigenvalues of the tangent map induced by its monodromy on the moduli space of the fundamental group of its fibre (see Main Theorem). 	This non abelian Reidemeister torsion is a combinatorial invariant of knots. 

\bigskip

	In~\cite{Fried:1988}, D. Fried has already computed the twisted Reidemeister torsion for bundles over the circle but in an \emph{acyclic} case. In the situation of fibered knots described in this article, we work in an \emph{non-acyclic} case. The key idea of our computations is to look at the Wang sequence in cohomology associated to the fibration. More precisely, we will compute the twisted Reidemeister torsion of the exterior of the fibered knot in terms of the Reidemeister torsion of the Wang sequence. In~\cite{LuckTorFib}, W. L\"uck, T. Schick and T. Thielmann study the behaviour of the analytic torsion under smooth fibrations. They obtain a general formula which involves several Reidemeister torsions, namely the torsions of the fibre,  of the basis and of the Leray-Serre spectral sequence for deRham cohomology induced by the fibration. This last term is of course the most difficult to compute. In our situation it coincides with the Wang sequence (see~\cite{Serre:1951}) and is precisely the one we will focus on and explicitly compute in terms of the monodromy. 
		
\bigskip
	
	The paper is organised as follows. Section~\ref{ReidemeisterTorsion}  reviews the sign-determined Rei\-de\-meis\-ter torsion. Section~\ref{ReidemeisterTorsionKnots} deals with the construction of the non abelian twisted Reidemeister torsion for knots. In Section~\ref{Fibered}, we prove the main theorem of the paper (see Main Theorem) about the twisted Reidemeister torsion associated to fibered knots. Finally, Section~\ref{Example} treats some examples.

\section{Review on the sign-determined Reidemeister torsion}
\label{ReidemeisterTorsion}

	The Reidemeister torsion of a finite simplicial complex $W$ is a more subtle invariant than the usual ones traditionally used in algebraic topology, because it uses an action of the fundamental group $\pi_1(W)$ on the universal cover of $W$. 
	This section reviews the basic definitions and sets up the conventions which will be used. For more details, we refer to Milnor's survey ~\cite{Milnor:1966} and to Turaev's monographs~\cite{Turaev:2000}~\&~\cite{Turaev:2002}. 

\subsection{Notation}
	In this paper, $\mathbb{F}$ is one of the fields $\IR$ or $\CC$; $\G$ is one of the Lie groups $\SU$ or $\SL$, and $\g$ is the associated Lie algebra $\su$ or $\sll$. Let $B_\g : \g \times \g \to \mathbb{F}$ denote the Killing form of $\g$. It is well-known that $B_\g$ is non-degenerate. 
	
	The Lie algebra $\su$ is identified with the \emph{pure quaternions}, i.e. with the quaternions of the form $q = a \ii + b \jj + c \kk$. In this case, $B_{\su}$ is equal to the usual scalar product $\langle \cdot, \cdot \rangle$ multiplied by $-2$. 
	
	The Lie algebra $\sll$ is identified with the space of $2 \times 2$ trace free matrices with complex entries. As a consequence, the Killing form satisfies
\[
B_{\sll}\left({\left(\begin{array}{cc}a & b \\c & -a\end{array}\right), \left(\begin{array}{cc}a' & b' \\c' & -a'\end{array}\right)}\right) = 8aa'+4(bc'+cb').
\]
	
\subsection{Basic definitions}
		Let $E$ be an $n$-dimensional vector space over $\mathbb{F}$. 
	For two ordered basis $\mathbf{a}  = \{ a_1, \ldots, a_n \}$ and $\mathbf{b} = \{ b_1, \ldots, b_n \}$ of $E$, we write $[\mathbf{a}/\mathbf{b}] = \det (p_{i j})_{i, j}$, where $a_i = \sum_{j=1}^n p_{i j} b_j$, for all $i$. 
	The bases $\mathbf{a}$ and $\mathbf{b}$ are called \emph{equivalent} if $[\mathbf{a}/\mathbf{b}]= +1$. In this case, we write $\mathbf{a} \sim \mathbf{b}$.
	
	Let $C_* = (\xymatrix@1@-.5pc{0 \ar[r] & C_n \ar[r]^-{d_n} & C_{n-1} \ar[r]^-{d_{n-1}} & \cdots \ar[r]^-{d_1} & C_0 \ar[r] & 0})$ be a chain complex of finite dimensional vector spaces over $\mathbb{F}$. For each $i$, consider $B_i = \im(d_{i+1} : C_{i+1} \to C_i)$, $Z_i = \ker(d_{i} : C_i \to C_{i-1})$ and the homology group $H_i = Z_i / B_i$. We say that $C_*$ is \emph{acyclic} if $H_i$ vanishes for all $i$. 

	 Suppose that for all $i = 1, \ldots, n$ both $C_i$ and $H_i$ are endowed with \emph{reference bases}. In this situation, $C_*$ is said to be \emph{based} and \emph{homology based}. One defines the Rei\-de\-meis\-ter torsion of such chain complexes as follows. Let $\mathbf{c}^i$ denote the reference basis for $C_i$ and let $\mathbf{h}^i$ denote the one for $H_i$. Let $\mathbf{b}^i$ be a sequence of vectors in $C_{i}$ such that $d_{i}(\mathbf{b}^i)$ is a basis of $B_{i-1}$ and let $\widetilde{\mathbf{h}}^i$ denote a lift of $\mathbf{h}^i$ in $Z_i$. For every $i$, the sequences of vectors $d_{i+1}(\mathbf{b}^{i+1})$, $\widetilde{\mathbf{h}}^i$ and $\mathbf{b}^{i}$ combines to yield a new basis $d_{i+1}(\mathbf{b}^{i+1})\widetilde{\mathbf{h}}^i\mathbf{b}^i$ of $C_i$. With this notation, the \emph{Reidemeister torsion} of $C_*$ (with reference bases $\mathbf{c}^*$ and $\mathbf{h}^*$) is the alternating product (see~\cite[Definition 3.1]{Turaev:2000}):
\begin{equation}
\label{Def:RTorsion}
\mathrm{tor}(C_*, \mathbf{c}^*, \mathbf{h}^*) = \prod_{i=0}^n [d_{i+1}(\mathbf{b}^{i+1})\widetilde{\mathbf{h}}^i\mathbf{b}^i/\mathbf{c}^i]^{(-1)^{i+1}} \in \mathbb{F}\setminus \{0\}.
\end{equation}
The torsion $\mathrm{tor}(C_*, \mathbf{c}^*, \mathbf{h}^*)$ does not depend on the choice of $\mathbf{b}^i$ and $\widetilde{\mathbf{h}}^i$. It does only depend on the equivalence classes of $\mathbf{c}^i$ and $\mathbf{h}^i$. More precisely, if $\mathbf{c'}^i$ is another basis of $C_i$ and $\mathbf{h'}^i$ another one of $H_i$, then we have the so-called \emph{basis change formula}
\begin{equation}
\label{EQ:changementdebase}
\frac{\mathrm{tor}(C_*, \mathbf{c'}^*, \mathbf{h'}^*)}{\mathrm{tor}(C_*, \mathbf{c}^*, \mathbf{h}^*)} = \prod_{i=0}^n \left( \frac{[{\mathbf{c}'}^i/\mathbf{c}^i]}{[{\mathbf{h}'}^i/\mathbf{h}^i]}\right)^{(-1)^i}.
\end{equation}

\subsection{The Reidemeister torsion of a CW-complex}

	If formula~(\ref{Def:RTorsion}) is used to define the Reidemeister torsion of a CW-complex, then we will fall into the well-known {up-to-sign ambiguity} of the Reidemeister torsion. To solve this problem, V.~Turaev has introduced a {sign-determined} Reidemeister torsion. 
	
\subsubsection*{\textbf{The sign-determined torsion}}

Keep the notation of the preceding subsection. Set 
\[
\alpha_i(C_*) = \sum_{k=0}^i \dim C_k \in \ZZ/2\ZZ, \; 
    \beta_i(C_*)  = \sum_{k=0}^i \dim H_k \in \ZZ/2\ZZ,
\]
\[
 |C_*| = \sum_{k\geqslant 0} \alpha_k(C_*) \beta_k(C_*) \in \ZZ/2\ZZ.
\]
The \emph{sign-determined Reidemeister torsion} of $C_*$ is the ``sign-corrected" torsion
\begin{equation}
\label{EQ:TorsionRaff}
	\mathrm{Tor}(C_*, \mathbf{c}^*, \mathbf{h}^*) = (-1)^{|C_*|}  \mathrm{tor}(C_*, \mathbf{c}^*, \mathbf{h}^*) \in \mathbb{F}\setminus \{0\},
\end{equation}
see~\cite[Section 3.1]{Turaev:1986} or~\cite[formula (1.a)]{Turaev:2002}. 
	
	For an acyclic based chain complex $C_*$, we have $\mathrm{Tor}(C_*) = \mathrm{tor}(C_*)$.

\subsubsection*{\textbf{The Reidemeister torsion of a CW-complex}}

	Let $W$ be a finite CW-complex; consider a representation $\rho \in \mathrm{Hom}(\pi_1(W); \G)$. The universal covering $\widetilde{W}$ of $W$ is endowed with the induced CW-structure and the fundamental group\,\footnote{We will not specify any base point, because all the constructions we do are invariant under conjugation.} $\pi_1(W)$ acts on $\widetilde{W}$ by the covering transformations. This action turns $C_*(\widetilde{W}; \ZZ)$ into a chain complex of left $\ZZ[\pi_1(W)]$-modules. The Lie algebra $\g$  can be viewed as a left $\ZZ[\pi_1(W)]$-module via the composition $Ad \circ \rho$, where $Ad : \SU \to \mathrm{Aut}(\su), A \mapsto Ad_A$ is the adjoint representation. Throughout the paper, we let $\g_{\rho}$ denote this $\ZZ[\pi_1(W)]$-module. 
The \emph{$\g_{\rho}$-twisted cochain complex} of $W$ is
\[
C^*(W; \g_\rho) = \mathrm{Hom}_{\pi_1(X)}(C_*(\widetilde{W}; \ZZ); \g_\rho).
\]
This cochain complex $C^*(W; \g_\rho)$ computes the \emph{$\g_\rho$-twisted cohomology} of $W$. We let $H^*(W; \g_\rho)$ denote the $\mathfrak{g}_\rho$-twisted cohomology of $W$. When $H^*(W; \g_\rho) = 0$, we say that $\rho$ is \emph{acyclic}.

	Choose a \emph{cohomology orientation} of $W$, \ie an orientation of the real vector space $H^*(W; \IR) = \bigoplus_{i\geqslant 0} H^i(W; \IR)$. Let $\mathfrak{o}$ denote this chosen orientation.
	
	 Let $\{e^{(i)}_1, \ldots, e^{(i)}_{n_i}\}$ denote the set of $i$-dimensional cells of $W$. Choose a lift $\tilde{e}^{(i)}_j$ of the cell $e^{(i)}_j$ in $\widetilde{W}$ and choose an arbitrary order and an arbitrary orientation for the cells $\tilde{e}^{(i)}_j$. Thus, for each $i$, $\mathbf{c}^{i} = \left\{ {\tilde{e}^{(i)}_1, \ldots, \tilde{e}^{(i)}_{n_i}} \right\}$ is a $\ZZ[\pi_1(W)]$-basis of $C_i(\widetilde{W}; \ZZ)$. If $\mathscr{B} = \{\mathbf{a}, \mathbf{b}, \mathbf{c}\}$ denotes an orthonormal basis of $\g$, then we consider the corresponding ``dual" basis 
	$$\mathbf{c}^{i}_{\mathscr{B}} = \left\{ \tilde{e}^{(i)}_{1, \mathbf{a}}, \tilde{e}^{(i)}_{1, \mathbf{b}}, \tilde{e}^{(i)}_{1, \mathbf{c}}, \ldots, \tilde{e}^{(i)}_{n_i, \mathbf{a}}, \tilde{e}^{(i)}_{n_i, \mathbf{b}}, \tilde{e}^{(i)}_{n_i, \mathbf{c}}\right\}$$ 
of $C^i(W; \mathfrak{g}_\rho) = \mathrm{Hom}_{\pi_1(X)}(C_*(\widetilde{W}; \ZZ); \g_\rho)$. If $\mathbf{h}^{i}$ denotes a basis for $H^i(W; \g_\rho)$, then $\mathrm{Tor}(C^*(W; \g_\rho), \mathbf{c}^*_{\mathscr{B}}, \mathbf{h}^{*}) \in \mathbb{F}\setminus \{0\}$ is well-defined.

	The cells $\{ \tilde{e}^{(i)}_j \}^{}_{0 \leqslant i \leqslant \dim W, 1 \leqslant j \leqslant n_i}$ are in one-to-one correspondence with the cells of $W$ and their order and orientation induce an order and an orientation for the cells $\{ e^{(i)}_j \}^{}_{0 \leqslant i \leqslant \dim W, 1 \leqslant j \leqslant n_i}$. We thus produce a basis over $\IR$ of $C^*(W; \IR)$ which is denoted by $c^*$. Provide each vector space $H^i(W; \IR)$ with a reference basis $h^i$ such that the basis $\left\{ {h^0, \ldots, h^{\dim W}} \right\}$ of $H^*(W; \IR)$ is {positively oriented} with respect to the cohomology orientation $\mathfrak{o}$. Compute the sign-determined Reidemeister torsion $\mathrm{Tor}(C^*(W; \IR), c^*, h^{*}) \in \mathbb{R}\setminus \{0\}$  of the resulting based and cohomology based chain complex and consider its sign $\tau_0 = \mathrm{sgn}\left(\mathrm{Tor}(C^*(W; \IR), c^*, h^{*})\right) \in \{\pm 1\}$.   
	The \emph{sign-determined Reidemeister torsion} of the cohomology oriented CW-complex $W$ is the product
\[
\mathrm{TOR}(W; \g_\rho, \mathbf{h}^{*}, \mathfrak{o}) = \tau_0 \cdot \mathrm{Tor}(C^*(W; \g_\rho), \mathbf{c}^*_{\mathscr{B}}, \mathbf{h}^{*}) \in \mathbb{F}\setminus \{0\}.
\]

	The torsion $\mathrm{TOR}(W; \g_\rho, \mathbf{h}^{*}, \mathfrak{o})$ is well-defined. It independent of the orthonormal basis $\mathscr{B}$ of $\g$. It does not depend  on the choice of the lifts $\tilde{e}^{(i)}_j$ nor on the choice of the chosen positively oriented basis of $H^*(W; \IR)$. Moreover, it is independent of order and orientation of the cells (because they appear twice). 
	Finally, it depends only on the conjugacy class of $\rho$.
	
	One can prove that $\mathrm{TOR}$ is invariant under cellular subdivision, homeomorphism class and simple homotopy equivalence. In fact, it is precisely the sign $(-1)^{|C_*|}$ in~(\ref{EQ:TorsionRaff}) which ensures all these important invariance properties hold (see~\cite[Chapter 2]{JDTHESE} for detailed proofs).
	
\subsection{The multiplicativity lemma}

	This lemma appears to be a very powerful tool for computing Reidemeister torsion. It will be used all over this paper.
\begin{MultLemma}\label{SS:Multiplicativite}
Let 
\begin{equation}\label{EQM}
\xymatrix@1@-.6pc{0 \ar[r] & C'_* \ar[r] & C_* \ar[r] & C_*'' \ar[r] & 0}
\end{equation}
be an exact sequence of chain complexes. Assume that $C'_*$, $C_*$ and $C''_*$ are based  and homology based. For all $i$, let ${\mathbf{c}'}^i$, ${\mathbf{c}}^i$ and ${\mathbf{c}''}^i$ denote the reference bases of $C_i'$, $C_i$ and $C_i''$ respectively. 
Associated  to~(\ref{EQM}) is the long sequence in homology 
\[
\xymatrix@1@-.6pc{ \cdots \ar[r] & H_i(C'_*) \ar[r] & H_i(C_*) \ar[r] & H_i(C''_*) \ar[r] & H_{i-1}(C'_*) \ar[r] & \cdots}
\] 
Let $\mathcal{H}_*$ denote this acyclic chain complex and base $\mathcal{H}_{3i+2} = H_i(C'_*)$, $\mathcal{H}_{3i+1} = H_i(C_*)$ and $\mathcal{H}_{3i} = H_i(C''_*)$ with the reference bases of  $H_i(C'_*)$, $H_i(C_*)$ and $H_i(C''_*)$ respectively.
If for all $i$, the bases ${\mathbf{c}'}^i$, ${\mathbf{c}}^i$ and ${\mathbf{c}''}^i$ are compatible, i.e. ${\mathbf{c}}^i \sim {\mathbf{c}'}^i{\mathbf{c}''}^i$, then
	$$\mathrm{Tor}(C_*) = (-1)^{\alpha(C'_*, C''_*) + \varepsilon(C'_*, C_*, C''_*)} \; \mathrm{Tor}(C'_*) \cdot \mathrm{Tor}(C''_*) \cdot \mathrm{tor}(\mathcal{H}_*),$$
where $$\alpha(C'_*, C''_*) = \sum_{i \geqslant 0} \alpha_{i-1}(C'_*) \alpha_i(C''_*) \in \ZZ/2\ZZ$$
and $$\varepsilon(C'_*, C_*, C''_*) = \sum_{i \geqslant 0} [(\beta_i(C_*)+1)(\beta_i(C_*')+\beta_i(C_*'')) + \beta_{i-1}(C_*')\beta_i(C_*'')] \in \ZZ/2\ZZ.$$
\end{MultLemma}
	The proof is a careful computation based on linear algebra, see~\cite[Lemma 3.4.2]{Turaev:1986} and~\cite[Theorem 3.2]{Milnor:1966}.

\section{Non abelian Reidemeister torsion for knot exteriors}
\label{ReidemeisterTorsionKnots}

	In this section, assume that $S^3$ is oriented and let $K \subset S^3$ be an oriented knot. Let $M_K = S^3 \setminus N(K)$ denote the exterior of $K$ and let $G_K = \pi_1(M_K)$ denote its group. Here $N(K)$ is an open tubular neighborhood of $K$. Recall that $M_K$ is a compact $3$-dimensional manifold with boundary $\bord M_K$ a $2$-dimensional torus. 
	The aim of this section is to construct a non abelian twisted Reidemeister torsion for $K$. For this purpose we shall produce some distinguished bases for the twisted cohomology group $H^*(M_K; \g_\rho)$. 
	
\subsection{Review on representation spaces}

	Given a finitely generated group $\pi$ we let $R(\pi; \G) = \mathrm{Hom}(\pi; \G)$ denote the space of $\G$-representation of $\pi$. This space is endowed with the {compact-open topology}. Here $\pi$ is assumed to have the discrete topology and the Lie group $\G$ is endowed with the usual one. A representation $\rho : \pi \to \G$ is called \emph{non-trivial} if $\rho(\pi) \not \subset \{\pm \I\}$. 
	
	A representation $\rho$ is called \emph{abelian} (resp. \emph{metabelian}) if its image $\rho(\pi)$ is an abelian subgroup of $\G$ (resp. if $\rho([\pi, \pi])$ is an abelian subgroup of $\G$). A representation $\rho : \pi \to \G$ is called \emph{reducible} if there exist a proper non-trivial subspace $U$ of $\CC^2$ such that $\rho(g)(U) \subset U$, for all $g \in \pi$. Observe that all abelian $\G$-representations of $\pi$ are reducible. A representation is called \emph{irreducible} if it is not reducible. We let $\widetilde{R}(\pi; \G)$ denote the subspace of irreducible representations. The Lie group $\G$ acts on $R(\pi; \G)$ by conjugation. We write $[\rho]$ for the conjugacy class of the representation $\rho \in R(\pi; \G)$. And we let  $\mathcal{R}(\pi; \G) = R(\pi; \G)/\G$ denote the moduli space.
	
	In the case of $\SU$-representations, observe that the reducible representations are exactly the abelian ones (which is not the case for $\SL$-representations). The action by conjugation of $\SU$ on $R(G; \SU)$ factors through $\SO=\SU/\{\pm \I\}$ as a \emph{free} action on the open  subspace $\widetilde{R}(\pi; \SU)$. Set $\widehat{R}(\pi; \SU) = \widetilde{R}(\pi; \SU)/\SO$. In this way, we can think of $\widehat{R}(\pi; \SU)$ as the base space of a principal $\SO$-bundle with total space $\widetilde{R}(\pi; \SU)$, see~\cite{GM:1992}, Section~3.A. 
	
	In the case of $\SL$-representations, the quotient $\mathcal{R}(\pi; \SL)$ is not Hausdorff in general. Following~\cite{CS:1983} we will focus on the \emph{representation variety} $X(\pi; \SL)$, which is the set of \emph{characters} of $\pi$. Associated to $\rho \in R(\pi, \SL)$ is its character $\chi_\rho : \pi \to \CC$, defined by $\chi_\rho(g) = \mathrm{Tr}(\rho(g))$. In some sense $X(\pi; \SL)$ is the ``algebraic quotient" of $R(\pi; \SL)$ by the action by conjugation of $\PSL$. We also let $\widehat{R}(\pi; \SL) = \widetilde{R}(\pi; \SL)/\PSL$ denote the image of $\widetilde{R}(\pi; \SL)$ under ${R}(\pi; \SL) \to X(\pi; \SL)$.
	
\subsection{Twisted cohomology of the torus}
\label{TwistT}
	
	Let $M$ be an $n$-dimensional compact manifold possibly with boundary $\bord M$. 
	If $M$ is oriented with an non-empty boundary, then $\bord M$ is a $(n-1)$-manifold which inherits an orientation by the convention ``\emph{the inward pointing normal vector in last position}". The Killing form $B_\g$ induces a \emph{cup-product}
\begin{equation}\label{cup}
\cup  \;: H^p(M; \g_\rho) \times H^{n - p}(M, \bord M; \g_\rho) \to H^{n}(M, \bord M; \mathbb{F}),
\end{equation}
which like $B_\g$ is  non-degenerate. 

	Let $T^2$ denote the $2$-dimensional torus. The following lemma computes the dimension of the $\g_\rho$-twisted cohomology groups of $T^2$.

\begin{lemma}\label{HT2}
Let $\rho \in R(\pi_1(T^2); \G)$ be a non-trivial. We have
\[
\dim_{\mathbb{F}} H^0(T^2; \g_\rho) = 1, \; \dim_{\mathbb{F}} H^1(T^2; \g_\rho) = 2 \text{ and } \dim_{\mathbb{F}} H^2(T^2; \g_\rho) = 1.
\] 
\end{lemma}
\begin{proof}	
	For any non-trivial representation $\rho \in R(\pi_1(T^2); \G)$ we observe that $$H^0(T^2; \g_\rho) = \g^{Ad \circ \rho(\pi_1(T^2))} \cong \mathbb{F}.$$ Poincaré duality implies $\dim_{\mathbb{F}} H^2(T^2; \g_\rho) = 1$. Finally, the fact that $$\sum_i (-1)^i \dim_{\mathbb{F}} H^i(T^2; \g_\rho) = 3\chi(T^2) = 0$$ furnishes $\dim_{\mathbb{F}} H^1(T^2; \g_\rho) = 2$.
\end{proof}
	
	 If $P^\rho$ denotes a generator of $H^0(T^2; \g_\rho)$, then the cup product~(\ref{cup}) and Poincaré duality combine to make the map $\phi^{(2)}_{P^\rho} : H^2(T^2; \g_\rho) \to H^2(T^2; \mathbb{F})$ given by $\phi^{2}_{P^\rho}(z) = P^\rho \cup z$. This is a natural isomorphism. 
	
	It is well-known that the non-trivial $\G$-representations of $\pi_1(T^2)$ are of two kinds: the \emph{hyperbolic} ones and the \emph{parabolic} ones. Here we say that the non-trivial representation $\rho$ of $\pi_1(T^2)$ in $\G$ is hyperbolic if each element in $\rho(\pi_1(T^2))$ is diagonalizable; $\rho$ is called parabolic if each element in $\rho(\pi_1(T^2))$ is non-diagonalizable or is $\I$. Further notice that each non-trivial element in ${R}(\pi_1(T^2); \SU)$ is hyperbolic.
	
	Assume that $\rho$ is hyperbolic. One can prove that the map $\phi^{(i)}_{P^\rho} : H^i(T^2; \g_\rho) \to H^i(T^2; \mathbb{F})$, defined by $\phi^*_{P^\rho}(z) = P^\rho \cup z$, is a natural isomorphism, for each $i=0, 1, 2$, because $B_\g(P^\rho, P^\rho) \ne 0$ (see~\cite[Proposition 3.18]{Porti:1997}). 
	If $\rho$ is assumed to be parabolic, then it can be shown that $\phi^{(1)}_{P^\rho} = P^\rho \cup \cdot$ is not an isomorphism.

\subsection{$\mu$-regular representations}
\label{mureg}

	We turn now to the case of knot groups. Let $\mu$ be a simple closed unoriented curve in $\bord M_K$. Among irreducible representations we focus on the $\mu$-regular ones. We say that $\rho \in \widetilde{R}(G_K; \G)$ is \emph{$\mu$-regular}, if (see~\cite[Definition 3.21]{Porti:1997}):
\begin{enumerate}
  \item the inclusion $\alpha : \mu \hookrightarrow M_K$ induces an \emph{injective} map $$\alpha^* : H^1(M_K; \g_\rho) \to H^1(\mu; \g_\rho),$$
  \item if $\mathrm{Tr}(\rho(\pi_1(\bord M_K))) \subset \{\pm 2\}$, then $\rho(\mu) \ne \pm \I$.
\end{enumerate}	
It is easy to see that this notion is invariant by conjugation. 

\begin{lemma}
If $\rho$ is $\mu$-regular, then $\dim_{\mathbb{F}} H^1(M_{K}; \g_{\rho}) = \dim_{\mathbb{F}} H^2(M_{K}; \g_{\rho}) = 1$.
\end{lemma}
\begin{proof}
	Assume that $\rho$ is $\mu$-regular. The inclusion $i : \bord M_K \hookrightarrow M_K$ induces an injective map $i^* : H^1(M_K; \g_\rho) \to H^1(\bord M_K; \g_\rho)$. We have $\dim_{\mathbb{F}} H^1(\bord M_K; \g_\rho) =2$ (see Lemma~\ref{HT2}) and Poincaré duality implies $\rk_{\mathbb{F}} i^* = 1$. Thus, $\dim_{\mathbb{F}} H^1(M_{K}; \g_{\rho}) = \rk_{\mathbb{F}} i^* = 1$. From $H^0(M_K; \g_\rho) = 0$ and $\sum_i (-1)^i \dim_{\mathbb{F}} H^i(M_K; \g_\rho) = 3\chi(M_K) = 0$ we deduce $\dim_{\mathbb{F}} H^2(M_K; \g_\rho) = 1$.
\end{proof}

	Here is an alternative formulation of $\mu$-regularity which will be more useful for us. The notion of $\mu$-regularity is defined for an unoriented curve $\mu$ but to avoid any ambiguity in what follows we must endow $\mu$ with a {``coherent" orientation}. 
	
	Let $\mathrm{int}(\cdot, \cdot)$ be the intersection form associated to the orientation of $\bord M_{K}$ induced by the one of $M_{K}$. The {peripheral subgroup} $\pi_1(\bord M_K)$ is generated by the meridian-longitude system $m, l$ of $K$. Here $m$ is oriented by the convention $\lk(K, m) = +1$ and $l$ is oriented by using the requirement that $\mathrm{int}(m, l) = +1$. We orient the curve $\mu$ as follows. If $\mu$ is parallel to $l$, then $\mu$ and $l$ are endowed with the {same} orientation, if not $\mu$ is endowed with the orientation such that $\mathrm{int}(\mu, l) \geqslant 0$. Let $\tilde{\mu}$ denote the resulting oriented curve.

Fix a generator $P^\rho$ of $H^0(\bord M_K; \g_\rho)$. The inclusion $\alpha : \mu \hookrightarrow M_K$ and the cup product~(\ref{cup}) induce the linear form $f^\rho_{\tilde{\mu}} : H^1(M_K; \g_\rho) \to \mathbb{F}$. We have $$f^\rho_{\tilde{\mu}}(v) = B_\g\left( {P^\rho, v(\tilde{\mu})}\right), \text{ for all } v \in  H^1(M_K; \g_\rho).$$ 
And we observe $f^\rho_{\tilde{\mu}^{-1}} = - f^\rho_{\tilde{\mu}}$, where $\tilde{\mu}^{-1}$ denotes the curve $\widetilde{\mu}$ with the opposite orientation. 

\begin{prop}
	The representation $\rho \in \widetilde{R}(G_K; \G)$ is $\mu$-regular if and only if the linear form $f^\rho_{\tilde{\mu}} : H^1(M_K; \g_\rho) \to \mathbb{F}$ is an isomorphism.
\end{prop}
\begin{proof}
If $\rho$ is $\mu$-regular, then $\rho_{|\pi_1(\bord M_K)}$ is non-trivial and we have $\dim_\mathbb{F} H^1(M_K; \g_\rho) = 1$. 
We split the first part of the proof into two steps to be clearer.
\begin{enumerate}
  \item If $\rho_{|\pi_1(\bord M_K)}$ is hyperbolic, then $H^*(\bord M_K; \g_\rho) \cong H^*(\bord M_K; \mathbb{F})$, thus the linear form $f^\rho_\mu$ is non-trivial.
  \item If $\rho_{|\pi_1(\bord M_K)}$ is parabolic then $\rho(\mu) \ne \pm \I$, thus $\g^{Ad \circ \rho(\mu)} = \g^{Ad \circ \rho(\pi_1(\bord M_K))}$. As a consequence $P^\rho \cup \cdot : H^1(\mu; \g_\rho) \to H^1(\mu; \mathbb{F})$ is an isomorphism, and thus $f^\rho_{\tilde{\mu}}$ is non-trivial.
\end{enumerate}

Assume now that $f^\rho_{\tilde{\mu}}$ is an isomorphism. We have $f^\rho_{\tilde{\mu}} = F \circ \alpha^*$, where $F : H^1(\mu; \g_\rho) \to \mathbb{F}$ is the linear form induced by the cup product. Thus $\alpha^*$ is injective, which proved the first assumption. Next, if $\rho(\mu) = \pm \I$, then $H^1(\mu; \g_\rho) \cong \g$ and thus $\rho_{|\pi_1(\bord M_K)}$ must be non-parabolic (see~\cite[Proposition 3.18]{Porti:1997}).
\end{proof}

\begin{example}\label{Ex}
	Let $K$ be a torus knot. All the irreducible representations of $G_K$ in $\G$ are $m$-regular and also $l$-regular (see~\cite[Example 1.43]{JDTHESE}).
\end{example}

\subsection{Reference bases of the twisted cohomology group of $M_K$}
\label{RefBasis2}
In this subsection assume that $\rho$ is $\mu$-regular. Fix a generator $P^\rho$ of $H^0(\bord M_K; \g_\rho)$ and keep the notation of Subsection~\ref{mureg}.

Firstly,  the \emph{reference generator} of $H^1(M_K; \g_\rho)$ is defined by
\begin{equation}\label{EQ:Defh1}
	h^{(1)}_\rho(\tilde{\mu}) = (f^\rho_{\tilde{\mu}})^{-1}(1).
\end{equation}

\begin{remark}
	The generator depends on the orientation of $\tilde{\mu}$, more precisely we have $h^{(1)}_\rho(\tilde{\mu}^{-1}) = - h^{(1)}_\rho(\tilde{\mu})$.
\end{remark}

Secondly, the construction of the reference generator of $H^2(M_K; \g_\rho)$ works as follows (see~\cite[Corollary 3.23]{Porti:1997}). The long exact sequence in $\g_\rho$-twisted cohomology associated to the pair $(M_K, \bord M_K)$ implies that the homomorphism $i^*: H^2(M_K; \g_\rho) \to H^2(\bord M_K; \g_\rho)$, induced by the inclusion $\bord M_K \hookrightarrow M_K$, is an isomorphism (because $\dim_{\mathbb{F}} H^2(M_K; \g_\rho) = \dim_{\mathbb{F}} H^2(\bord M_K; \g_\rho)  = 1$). As a consequence, the composition $$\phi^{(2)}_{P^\rho} \circ i^* : H^2(M_K; \g_\rho) \to H^2(\bord M_K; \g_\rho) \to H^2(\bord M_K; \mathbb{F}) = H^2(\bord M_K; \ZZ) \otimes \mathbb{F}$$ is an {isomorphism}. Let $c$ be the generator of $H^2(\bord M_K; \ZZ) = \mathrm{Hom}(H_2(\bord M_K; \ZZ), \ZZ)$ 
corresponding to the fundamental class $\lbrack \! \lbrack \bord M_K \rbrack \! \rbrack \in H_2(\bord M_K; \ZZ)$ induced by the orientation of $\bord M_K$. The \emph{reference generator} of $H^2(M_K; \g_\rho)$ is defined by 
\begin{equation}\label{EQ:Defh2}
h^{(2)}_\rho = (\phi^{(2)}_{P^\rho} \circ i^*)^{-1}(c).
\end{equation}

\subsection{The Reidemeister torsion $\mathbb{T}^K_\mu$}
\label{defTau}

	We equip the exterior of $K$ with its \emph{canonical cohomology orientation} defined as follows (see~\cite[Section V.3]{Turaev:2002}). We have 
	$$H^*(M_K; \IR) = H^0(M_K; \IR) \oplus H^1(M_K; \IR)$$ 
and we base this $\IR$-vector space with $\{ \lbrack \! \lbrack pt \rbrack \! \rbrack, m^*\}$. Here $\lbrack \! \lbrack pt \rbrack \! \rbrack$ is the cohomology class of a point, and $m^* : m \mapsto 1$ is the dual of the meridian $m$ of $K$. This reference basis of $H^*(M_K; \IR)$ induces the so-called canonical cohomology orientation of $M_K$. In the sequel, we let $\mathfrak{o}$ denote the canonical cohomology orientation of $M_K$.
	
	Let $\rho : G_K \to \G$ be a $\mu$-regular representation. The \emph{Reidemeister torsion $\mathbb{T}^K_\mu$} at $\rho$ is  
defined by 
\[
\mathbb{T}^K_\mu(\rho)  = \mathrm{TOR}\left( {M_K; \g_\rho, \{h^{(1)}_\rho(\tilde{\mu}), h^{(2)}_\rho\}, \mathfrak{o}} \right).
\]

	Observe that the $\mu$-torsion $\mathbb{T}^K_\mu(\rho)$ does not depend on the choice of the generator $P^\rho$ of $H^0(\bord M_K; \g_\rho)$. This property is a consequence of formula~(\ref{EQ:changementdebase}), because $h^{(1)}_\rho(\tilde{\mu})$ and $h^{(2)}_\rho$ change in the same way at the same time when we change $P^\rho$.

\begin{prop}
	The torsion $\mathbb{T}^K_\mu(\rho)$ does not depend on the orientation of $K$. 
\end{prop}
\begin{proof}
If we change the orientation of $K$, then the orientations of $m$ and $l$ change simultaneously; thus the orientation of $\tilde{\mu}$ is reversed. As a consequence, the reference generator of $H^2(M_K; \g_\rho)$ is unchanged, but the one of $H^1(M_K; \g_\rho)$ and the cohomology orientation are reversed simultaneously. Thus $\mathbb{T}^K_\mu(\rho)$ does not change.
\end{proof}

\begin{remark}
For a hyperbolic knot $K$, the torsion $\mathbb{T}^K_\mu(\rho)$ is a sign-refined version of the inverse of Porti's torsion function (see~\cite{Porti:1997}).
\end{remark}

\section{The Reidemeister torsion for fibered knots}
\label{Fibered}

	In this section, we give an explicit formula to compute the sign-determined non abelian Reidemeister torsion associated to fibered knots using the monodromy of the knot. 
	
	Let $K \subset S^3$ be a fibered knot. We let $F$ denote the fiber of $K$ and let $\gamma$ denote the boundary of $F$. Recall that $\gamma$ corresponds to the longitude of $K$ and thus is an oriented curve in $\bord M_K$ (see Subsection~\ref{mureg}). If $g$ denotes the genus of the surface $F$, then $\widehat{R}(\pi_1(F); \G)$ is smooth and $\dim_{\mathbb{F}}\widehat{R}(\pi_1(F); \G) = 6g-3$. Let $\varphi = \rho_{|\pi_1(F)}$ be the restriction of $\rho$ to $\pi_1(F)$. The monodromy $\phi : F \to F$ induces a diffeomorphism $R(\phi) : \widehat{R}(F; \G) \to \widehat{R}(F; \G)$. Let $I_\gamma : \widehat{R}(F; \G) \to \mathbb{F}$, be the function defined by $I_\gamma : \varrho \mapsto \mathrm{Tr}(\varrho(\gamma))$. As $\phi$ preserves the boundary $\gamma = \bord F$ we see that $I_\gamma \circ R(\phi) = I_\gamma$, because the trace function is invariant by conjugation.		

\begin{maintheorem*}
	Let $K \subset S^3$ be a fibered knot. Let $F$ be its fiber, a surface of genus $g$ and boundary $\gamma$. Assume that $\varepsilon_0$ is the sign of the determinant of the isomorphism $\mathrm{Id} - \phi^* : H^1(F; \IR) \to H^1(F; \IR)$, where $\phi^*$ is induced by the monodromy. If $\rho : G_K \to \G$ is a non-metabelian $\gamma$-regular representation, then the tangent map at $\varphi = \rho_{|\pi_1(F)}$ to $R(\phi) : \widehat{R}(F; \G) \to \widehat{R}(F; \G)$ admits $1$ as simple eigenvalue. If we let $\lambda_1, \ldots, \lambda_{6g-4}$ denote its other eigenvalues, then
\begin{equation}
\label{E:TorsionFibre}
\mathbb{T}^K_\gamma(\rho) =  - \varepsilon_0 \cdot \prod_{i=1}^{6g-4} \frac{1}{1-\lambda_i} \in \mathbb{F}\setminus \{0\}.
\end{equation}
\label{T:TorsionFibre}
\end{maintheorem*}

	In the Main Theorem, we restrict our attention on the irreducible non-metabelian representations of $G_K$ in $\G$. This restriction is just technical; moreover, there exists only a finite number of irreducible metabelian representations of $G_K$ in $\SU$ (see~\cite[Proposition 4.2]{Lin:2001}).

	The main tool to compute the $\gamma$-torsion associated to the fibered knot $K$ with coefficients in $\g_\rho$ is the Wang sequence in twisted cohomology associated to the fibration $F \hookrightarrow M_K \to S^1$. This idea is similar to the one used by D. Fried in~\cite{Fried:1988} in an acyclic case but it is technically different, because of the non triviality of $H^*(M_K; \g_\rho)$. 	

	Let $a_1, b_1, \ldots, a_g, b_g$ denote the generators of $\pi_1(F)$. If $\phi_* : \pi_1(F) \to \pi_1(F)$ is the homomorphism induced by the monodromy $\phi$, then $G_K$ admits the presentation
\begin{equation}\label{G_fibre}
G_K = \langle a_1, b_1, \ldots, a_g, b_g, t \; |\; t^{-1} a_i t= \phi_*(a_i), t^{-1} b_i t = \phi_*(b_i), \; 1 \leqslant i\leqslant g \rangle
\end{equation}
where $t$ corresponds to the meridian of $K$. We know that $M_K$ collapses to a $2$-dimensional CW-complex. 
More precisely, the presentation~(\ref{G_fibre}) of $G_K$ allows us to define a $2$-dimensional CW-complex $X_K$ as follows. The $0$-skeleton of $X_K$ consists of one point, the $1$-skeleton $X_K^1$ is a wedge of $2g+1$ oriented circles corresponding to the generators of~(\ref{G_fibre}); finally $X_K$ is obtained from $X^1_K$ by gluing $2g$ closed $2$-cells attached using the relations of~(\ref{G_fibre}). A result of Waldhausen~\cite{Waldhausen:1978} implies that the cell complexes $X_K$ and $M_K$ have the same simple homotopy type. As a consequence $X_{K}$ will be used to explicitly compute $\mathbb{T}^K_\gamma(\rho)$.         Remark that $H^*(X_K; \g_\rho)  = H^*(M_K; \g_\rho)$.

\section{Proof of the Main Theorem}\label{Proof}

With the notation of the previous section, we have $$\mathbb{T}^K_\gamma(\rho) = \tau_0 \cdot \mathrm{Tor}(X_K; \g_\rho, \{h^{(1)}_\rho(\gamma), h^{(2)}_\rho\}),$$
where $\tau_0$ is the sign of $\mathrm{Tor}(X_K; \IR)$. The proof of  the Main Theorem consists in the computation of the torsions $\mathrm{Tor}(X_K; \g_\rho)$ and $\mathrm{Tor}(X_K; \IR)$ in terms of the torsions of the Wang sequences with twisted and with real coefficients respectively. It is divided into several steps.

\subsection{Preliminaries: The Wang sequence with twisted coefficients}
\label{WangS}

	The monodromy $\phi : F \to F$ induces an action $\phi^*_\rho : C^*(F; \g_\varphi) \to C^*(F; \g_\varphi)$ at the level of the twisted chain complex of $F$. 
	Thus  
\begin{equation}\label{action}
\xymatrix@1{0 \ar[r] & C^*(M_K; \g_\rho) \ar[r]^-{i^*} & C^*(F; \g_\varphi) \ar[r]^-{\mathrm{Id}-\phi^*_\rho} & C^*(F; \g_\varphi) \ar[r] & 0,}
\end{equation}
is an exact sequence of chain complexes.

	Observe that the representation $\varphi$ is non-abelian, because $[G_K, G_{K}] = \pi_1(F)$ and $\rho$ is supposed to be non-metabelian. Thus, $H^0(F; \g_\varphi) = \g^{Ad \circ \varphi(\pi_1(F))} = 0$ (see~\cite[Lemma 0.7]{Porti:1997}). As a consequence sequence~(\ref{action}) induces the following long exact sequence in twisted cohomology: 
\begin{equation}\label{Wang}
\mathscr{W}^\rho_* = \xymatrix@1{0 \ar[r] & H^1(M_K; \g_\rho) \ar[r]^-{i^*} & H^1(F; \g_{\varphi}) \ar[r]^-{\mathrm{Id} - \phi^*_\rho} & H^1(F; \g_{\varphi}) \ar[r]^-\delta & H^2(M_K; \g_\rho) \ar[r] & 0.}
\end{equation}
Sequence~\ref{Wang} is called the \emph{Wang sequence} (with twisted coefficients) associated to the fibration $F \hookrightarrow M_K \to S^1$.

	In sequence~(\ref{Wang}), observe that $\phi^*_\rho : H^1(F; \g_\varphi) \to H^1(F; \g_\varphi)$ can be identified with the tangent map at $\varphi$ to $R(\phi)$, because $H^1(F; \mathfrak{g}_\varphi \cong \tangent{\varphi}{\widehat{R}(\pi_1(F); \G)}$ (see~\cite[\S~3.1.3.]{Porti:1997} and~\cite[\S~1.4.4. p.~17]{JDTHESE}). Explicitly, we have $\phi^*_\rho(h)(x) = Ad_{\rho(m)^{-1}}(h(\phi_*(x)))$, for all $h \in H^1(F; \mathfrak{g}_\varphi)$ and for all $x \in \pi_1(F)$. 
	The equality $I_\gamma \circ R(\phi) = I_\gamma$ implies that $1$ is always an eigenvalue of $\phi^*_\rho$. Combining the exactness of the Wang sequence~(\ref{Wang}) to the fact that $\rho$ is $\gamma$-regular, we further observe that $1$ is a simple eigenvalue because $\dim_{\mathbb{F}} \ker(\mathrm{Id} - \phi^*_\rho) = \dim_{\mathbb{F}} H^1(M_K; \mathfrak{g}_\rho) = 1$. 

We compute the $\mathfrak{g}_\rho$-twisted Reidemeister torsion $\mathrm{Tor}(X_K; \g_\rho)$ as follows. The twisted cohomology group $H^1(M_K; \g_\rho) \cong \mathbb{F}$ is based with the generator $h^{(1)}_\rho(\gamma)$ (cf.~equation~(\ref{EQ:Defh1})) and the twisted cohomology group $H^2(M_K; \g_\rho) \cong \mathbb{F}$ is based with the generator $h^{(2)}_\rho$ (cf.~equation~(\ref{EQ:Defh2})). Next, fix a chosen basis (over $\mathbb{F}$) of $H^1(F; \g_\varphi)$ and observe that $\mathrm{Tor}(\mathscr{W}^\rho_*)$ does not depend on this choice (see equation~(\ref{EQ:changementdebase})). Further notice that it is precisely this indeterminacy which will be used to compute $\mathrm{Tor}(X_K; \g_\rho)$ in terms of $\mathrm{Tor}(\mathscr{W}^\rho_*)$.
We have
$$\alpha = \alpha(C^*(X_K; \g_\rho), C^*(F; \g_\varphi))= 1 \in \ZZ/2\ZZ$$ and $$\varepsilon =\varepsilon(C^*(X_K; \g_\rho), C^*(F; \g_\varphi), C^*(F; \g_\varphi))= 0 \in \ZZ/2\ZZ.$$ 
As a consequence, the multiplicativity lemma applied to sequence~(\ref{action}) provides
\begin{equation}\label{TorsionWang}
\mathrm{Tor}(X_K; \g_\rho, \{h^{(1)}_\rho(\gamma), h^{(2)}_\rho\}) = - (\mathrm{Tor}(\mathscr{W}^\rho_*))^{-1}.
\end{equation}

\subsection{Torsion of the Wang sequence with twisted coefficients}
The aim of this subsection is to show
\begin{claim}
Under the hypothesis of the Main Theorem, we have
\begin{equation}\label{E:Torsion2}
\mathrm{Tor}(\mathscr{W}^\rho_*) = \prod_{i=1}^{6g-4} (1 - \lambda_i) \in \mathbb{F}\setminus \{0\}.
\end{equation}
\end{claim}
\begin{proof}
	Recall that $\rk_{\mathbb{F}}(\mathrm{Id} - \phi^*_\rho) = 6g-4$. Choose a basis $\mathbf{c} = \left\{{v_1, v_2, \ldots, v_{6g-3}}\right\}$  for $H^1_\varphi(F)$ in which $\phi^*_\rho$ is upper-triangular and such that $i^*(h^{(1)}_\rho(\gamma)) = v_{6g-3}$. In the sequel, we assume that $H^1(F; \g_{\varphi})$ is based with $\mathbf{c}$. 

	Let $\lambda_1, \ldots, \lambda_{6g-4}$ denote the eigenvalues of $\phi^*_\rho$ different from $1$. With this notation, we have:
\begin{enumerate}
  \item $\phi^*_\rho(v_{6g-3}) = v_{6g-3}$ and $v_{6g-3} \notin \im(\mathrm{Id} - \phi^*_\rho)$,
  \item $ \left\{{(\mathrm{Id} - \phi^*_\rho)(v_1), \ldots, (\mathrm{Id} - \phi^*_\rho)(v_{6g-4})}\right\}$ is a basis for $\im(\mathrm{Id} - \phi^*_\rho)$,
  \item $\delta(v_{6g-3}) = h^{(2)}_\rho$. 
\end{enumerate}

  	The first two statement are easy. The third one is obtained as follows. Combine the Wang sequences (with twisted coefficients) associated to the fibration $F \hookrightarrow M_K \to S^1$ with the one associated to the fibration $\gamma \hookrightarrow \bord M_K \to S^1$ to obtain the commutative diagram: 
 \[
\xymatrix@-.7pc{0 \ar[r] &H^1(M_K; \g_\rho) \ar[r] \ar[d] & H^1(F; \g_{\varphi}) \ar[r] \ar[d] & H^1(F; \g_{\varphi}) \ar[r]^-\delta \ar[d] & H^2(M_K; \g_\rho) \ar[d]^-\cong \ar[r] & 0 \\ 
\cdots \ar[r] & H^1(\bord M_K; \g_\rho) \ar[r]  & H^1(\gamma; \g_{\varphi}) \ar[r] & H^1(\gamma; \g_{\varphi}) \ar[r]^-\Delta_-\cong & H^2(\bord M_K; \g_\rho) \ar[r] & 0} 
\]
In the preceding diagram, the vertical arrows are induced by inclusions.

If $\mathbf{b}^1 = \left\{{v_{6g-3}}\right\}$ and $\mathbf{b}^2 = \left\{v_1, \ldots, v_{6g-4}\right\}$, then
\[\mathrm{Tor}(\mathscr{W}^\rho_*) = \left[{(\mathrm{Id} - \phi^*_\rho)(\mathbf{b}^2)\mathbf{b}^1 / \mathbf{c}}\right] = \prod_{i=1}^{6g-4} (1 - \lambda_i).\]
\end{proof}

\subsection{Torsion of the Wang sequence with real coefficients}

	For the same reason as in Subsection~\ref{WangS}, the fibration $F \hookrightarrow M_K \to S^1$ induces a Wang sequence $\mathfrak{W}_*$ with real coefficients which splits into three isomorphisms:
\begin{equation}\label{iso}
\xymatrix@1@-.9pc{H^0(M_K; \IR) \ar[r]^-{i^*}_-{\cong} & H^0(F; \IR)}, \; \xymatrix@1@-.9pc{H^0(F; \IR) \ar[r]^-{\delta}_-{\cong} & H^1(M_K; \IR)} \text{ and } \xymatrix@1{H^1(F; \IR) \ar[r]^-{\mathrm{Id - \phi^*}}_-{\cong} & H^1(F; \IR).}
\end{equation}

	Observe that $\mathrm{Tor}(\mathfrak{W}_*)$ does not depend on the choice of the bases of $H^*(F; \IR)$ used for the computation (see equation~(\ref{EQ:changementdebase})). As a consequence, we will choose appropriate bases to compute the torsion. Recall that $H^0(M_K; \IR)$ is based with the generator $\lbrack \! \lbrack pt \rbrack \! \rbrack$ and $H^1(M_K; \IR)$ is based with $m^*$ (see Subsection~\ref{defTau}). Suppose that $H^0(F; \IR)$ is endowed with the generator $i^*(\lbrack \! \lbrack pt \rbrack \! \rbrack)$. We fix an arbitrarily chosen basis of $H^1(F; \IR)$. 
	
We have
$$\alpha = \alpha(C^*(X_K; \IR), C^*(F; \IR)) = 1 \in \ZZ/2\ZZ$$
and
$$\varepsilon =\varepsilon(C^*(X_K; \IR), C^*(F; \IR), C^*(F; \IR)) =1 \in \ZZ/2\ZZ.$$
The multiplicativity lemma provides
\[
\mathrm{Tor}(X_K; \IR) =  (\mathrm{Tor}(\mathfrak{W}_*))^{-1},
\]
so that
\begin{equation}\label{Tau}
\tau_0 = \mathrm{sgn}(\mathrm{Tor}(X_K; \IR)) = \mathrm{sgn}(\mathrm{Tor}(\mathfrak{W}_*)).
\end{equation}
	
	We turn now to the computation of $\mathrm{Tor}(\mathfrak{W}_*)$. More precisely we show
	
\begin{claim} We have
\begin{equation}\label{TorR}
\mathrm{sgn}(\mathrm{Tor}(\mathfrak{W}_*)) = \varepsilon_0.
\end{equation}
\end{claim}
Here we recall that $\varepsilon_0 = \mathrm{sgn}(\det(\mathrm{Id} - \phi^*))$, where $\phi^* : H^1(F; \IR) \to H^1(F; \IR)$ is induced by the monodromy.
\begin{proof}
	Notice that $\mathrm{Tor}(\mathfrak{W}_*)$ is the product of the determinants of the three isomorphisms~(\ref{iso}) in the bases chosen above.
	
	First, in the chosen bases, the determinant of the isomorphism $i^* : H^0(X_K; \IR) \to H^0(F; \IR)$ is $1$.
	
	Second, we prove that $$\mathrm{sgn}([\delta \circ i^*(\lbrack \! \lbrack pt \rbrack \! \rbrack) / m^*])  = 1.$$

	In fact, we show that $\delta(i^*(\lbrack \! \lbrack pt \rbrack \! \rbrack)) = m^*$.
	Combine the Wang sequences in cohomology (with coefficients in $\ZZ$) associated to the fibrations $F \hookrightarrow M_K \to S^1$ and $\gamma \hookrightarrow \bord M_K \to S^1$ to obtain the commutative diagram
\[
\xymatrix@-.9pc{0 \ar[r] & H^0(M_K; \ZZ) \ar[r]^-{i^*} \ar[d]^-\cong & H^0(F; \ZZ) \ar[r]^-0 \ar[d]^-\cong & H^0(F; \ZZ) \ar[r]^-\delta \ar[d]^-\cong & H^1(M_K; \ZZ) \ar[r]^-{i^*} \ar[d]^-{j^*} & H^1(F; \ZZ) \ar[r] \ar[d] & \cdots \\
0 \ar[r] & H^0(\bord M_K; \ZZ) \ar[r]^-\cong & H^0(\gamma; \ZZ) \ar[r]^-0 & H^0(\gamma; \ZZ) \ar[r]^-\Delta & H^1(\bord M_K; \ZZ) \ar[r]^-{\iota^*} & H^1(\gamma; \ZZ) \ar[r]^-0 &\cdots}
\]
In the preceding diagram, the vertical arrows are induced by inclusions.

The group $H^0(\bord M_K; \ZZ)$ is generated by $\lbrack \! \lbrack pt \rbrack \! \rbrack$, the group $H^0(\gamma; \ZZ)$ by $i^*(\lbrack \! \lbrack pt \rbrack \! \rbrack)$ and $H^1(\bord M_K; \ZZ)$ is based with $\left\{{m^*, \gamma^*}\right\}$; here $m$ denotes the meridian of $K$ and $\gamma = \bord F$ its longitude. We know that $j^*(m^*) = m^*$.

Using the previous commutative diagram, to prove $\delta(i^*(\lbrack \! \lbrack pt \rbrack \! \rbrack)) = m^*$ it is enough to show $\Delta(i^*(\lbrack \! \lbrack pt \rbrack \! \rbrack)) = m^*$. This last equality is obtained by a careful examination of the short exact sequence
\begin{equation}\label{gamm}
\xymatrix@1@-.7pc{0 \ar[r] & H^0(\gamma; \ZZ) \ar[r]^-\Delta & H^1(\bord M_K; \ZZ) \ar[r]^-{\iota^*} & H^1(\gamma; \ZZ) \ar[r] & 0.}
\end{equation}
In sequence~(\ref{gamm}), $H^0(\gamma; \ZZ)$ is generated by $i^*(\lbrack \! \lbrack pt \rbrack \! \rbrack)$, $H^1(\gamma; \ZZ)=\mathrm{Hom}(H_1(\gamma; \ZZ), \ZZ)$ by $\gamma^* : \gamma \mapsto 1$ (which is the dual of $\gamma$) and $H^1(\bord M_K; \ZZ) = \mathrm{Hom}(H_1(\bord M_K; \ZZ), \ZZ)$ is based with $\left\{{m^*, \gamma^*}\right\}$, where $\iota^*(\gamma^*) = \gamma^*$.
\end{proof}

\subsection{End of the proof}
It remains to bring together all the computations we have done before.

\begin{proof}[Proof of the Main Theorem]
	Recall that $$\mathbb{T}^K_\gamma(\rho) = \tau_0 \cdot \mathrm{Tor}(X_K; \g_\rho, \{h^{(1)}_\rho(\gamma), h^{(2)}_\rho\}).$$  
	
	On the one hand, equations~(\ref{TorsionWang}) and~(\ref{E:Torsion2}) imply $$\mathrm{Tor}(X_K; \g_\rho, \{h^{(1)}_\rho(\gamma), h^{(2)}_\rho\}) = - \prod_{i=1}^{6g-4} \frac{1}{1-\lambda_i}.$$
	 
	 On the other hand, equations~(\ref{Tau}) and~(\ref{TorR}) imply $\tau_0 = \varepsilon_0$, and this completes the proof.
\end{proof}

\section{Examples} 
\label{Example}

This last section is devoted to concrete computations. The Main Theorem can be applied to provide explicit formulae. We focus our attention first on  the $\SU$-representation space of the group of torus knots and next on the $\SL$-repre\-sen\-ta\-tion space of the group of the figure eight knot. 

\subsection{The trefoil knot}
\label{trefle}
	It is well-known that the (right hand) trefoil knot $K$ is a fibered knot of genus $1$. Let $F$ denote the fiber of $K$ and $\gamma = \bord F$ its longitude. We know that $F$ is a one-punctured $2$-dimensional surface of genus $1$. In the sequel on this subsection, we let $a, b$ denote the generators of the free group $\pi_1(F)$. The group $G_K$ of $K$ admits, as a fibered knot, the following presentation: 
\[
G_K = \langle a, b, t \; |\; t^{-1}at = ab^{-1}a^{-1}, t^{-1}bt = ab \rangle,
\]
in which $t$ represents the meridian of $K$. 

	Recall that $\widehat{R}(G_K; \SU)$ is the set of $\SO$-conjugacy classes of the irreducible representations $\rho : G_K \to \SU$ such that $\mathrm{Tr}(\rho(t)) = \sqrt{3} \cos(\theta)$, for $\theta \in (0, \pi)$, see~\cite[Theorem 1]{Klassen:1991}. Observe  that each irreducible representations of $G_K$ in $\SU$ is $\gamma$-regular (see Example~\ref{Ex}).

	Recall that $H_1(F; \ZZ) = [\pi_1(F), \pi_{1}(F)]$; $H^1(F; \ZZ) = \mathrm{Hom}(H_1(F; \ZZ), \ZZ)$ is endowed with the basis $\{a^*, b^*\}$, where $a^* : a \mapsto 1$, $b^* : b \mapsto 1$ respectively denote the dual of $a$ and the one of $b$. The monodromy $\phi : F \to F$ induces an endomorphism $\phi^* : H^1(F; \ZZ) \to H^1(F; \ZZ)$ such that the matrix of $\phi^*$ in the basis $\{a^*, b^*\}$ is $\left(\begin{array}{cc}0 & 1 \\-1 & 1\end{array}\right)$. As a consequence, $\varepsilon_0 = \det (\mathrm{Id} - \phi^*) = +1$.

Let $\rho$ be an irreducible $\SU$-representation of $G_K$. 
Set $x_1 = I_a$, $x_2 = I_b$ and $x_3 = I_{ab}$. We know that the moduli space $\widehat{R}(F; \SU)$ is parametrized by $x_1, x_2, x_3$. Observe that, with respect to the coordinates $(x_1, x_2, x_3)$, the action of $\phi^*_\rho$ is given by 
	$$P = (x_2, x_3, x_1) \in \ZZ[x_1, x_2, x_3]^3.$$ 
Thus the tangent map at $\varphi = \rho_{|\pi_1(F)}$ to $R(\phi)$ is given by the matrix
\begin{equation}\label{ATref}
 \left( \frac{\bord P_i}{\bord x_j}(\varphi)\right)_{i, j} = \left(\begin{array}{ccc}0 & 0 & 1 \\1 & 0 & 0 \\0 & 1 & 0\end{array}\right).
 \end{equation}
The matrix~(\ref{ATref}) admits $1$, $e^{2i\pi/3}$ and $e^{-2i\pi/3}$ as eigenvalues.
As a consequence, the Main Theorem gives 
	$$\mathbb{T}^K_\gamma(\rho) = -\frac{1}{3}.$$

\subsection{The torus knots}
	More generally, let $p, q \in \mathbb{N}^*$ be coprime and let $K_{p, q}$ denote the (right hand) torus knot of type $(p, q)$. Its group admits the presentation $G_{p, q} = \langle x, y \; |\; x^p = y^q \rangle$. It is known that $K_{p, q}$ is fibered of genus $(p-1)(q-1)/2$ (see~\cite{BZ:1996}). 
	
	Let $\rho : G_{p, q} \to \SU$ be an irreducible representation such that $\mathrm{Tr}(\rho(x)) = 2\cos\left(\frac{\pi a}{p}\right)$,  $\mathrm{Tr}(\rho(y)) = 2\cos\left(\frac{\pi b}{q}\right)$, where $0 < a < p$, $0< b< q$ and $a \equiv b \,(\mathrm{mod.} 2)$. It is well-known that $\widehat{R}(G_{p, q}; \SU)$ consists of the conjugacy classes of all such representations (see~\cite[Theorem 1]{Klassen:1991}). In a similar way as for the trefoil knot, $\rho$ is $\gamma$-regular and we obtain
$$\mathbb{T}^{K_{p,q}}_{\gamma}(\rho) = - \frac{16}{p^2q^2} \sin^2\left(\frac{\pi a}{p}\right) \sin^2\left(\frac{\pi b}{q}\right).$$

\subsection{The figure eight knot}
	We turn now to the case of a hyperbolic knot: the figure eight knot. We let $K$ denote it and let $F$ denote its fiber. As a fibered knot, the group $G_K$ of the figure eight knot $K$ admits the following presentation:
\[
G_K = \langle a, b, t \; | \; t^{-1} a t = ab, t^{-1}bt = bab \rangle.
\]	
We study the moduli space of $\SL$-representations of $G_K$. 
The monodromy $\phi : F \to F$ induces an endomorphism $\phi^* : H^1(F; \ZZ) \to H^1(F; \ZZ)$ such that the matrix of $\phi^*$ in the basis $\{a^*, b^*\}$ is $\left(\begin{array}{cc}1 & 1 \\ 1 & 2\end{array}\right)$. As a consequence, $\varepsilon_0 = \det(\mathrm{Id} - \phi^*) = -1$. With the same notation as in Subsection~\ref{trefle}, the action of $\phi^*_\rho$ is given by 
$$P = (x_3, x_2x_3 - x_1, x_2x_3^2-x_1x_3-x_2) \in \ZZ[x_1, x_2, x_3]^3.$$
Further notice that for $\rho \in \widetilde{R}(G_K; \SL)$ we have 
\begin{equation}\label{eq}
I_a(\rho) = I_c(\rho) \text{ and } I_a(\rho) + I_b(\rho) = I_a(\rho)I_b(\rho).
\end{equation}
 Combining the Main Theorem to equalities~(\ref{eq}) we obtain (compare with~\cite[p. 113]{Porti:1997})	
$$\mathbb{T}^K_\gamma(\rho) = \frac{ 1}{3-2(I_a(\rho)+I_b(\rho))}.$$
The well-known identity $$\mathrm{Tr}(ABA^{-1}B^{-1}) = -2 - \mathrm{Tr}(A)\mathrm{Tr}(B)\mathrm{Tr}(AB) + (\mathrm{Tr}(A))^2 +(\mathrm{Tr}(B))^2+(\mathrm{Tr}(AB))^2$$ implies $I_\gamma = x_1^2+x^2_2-x_1-x_2-2$. As a consequence,
$${\mathbb{T}^K_\gamma(\rho)}^2 = \frac{1}{{17+4 I_\gamma(\rho)}}.$$

The hyperbolic structure of the exterior of $K$ determines a unique (up to complex-conjugation) \emph{holonomy representation} in $\mathrm{PSL}_2(\CC)$ which lifts to two representations $\varrho_{\pm}$ in $\SL$ satisfying
\[
\varrho_{\pm}(m) = \left(\begin{array}{cc}\pm 1 & 1 \\0 & \pm 1\end{array}\right) \text{ and } \varrho_{\pm}(l) = \left(\begin{array}{cc}1 & \pm 2i\sqrt{3} \\0 & 1\end{array}\right).
\] 
Thus $$\mathbb{T}^K_\gamma(\varrho_{\pm}) = \frac{1}{5}.$$

\section*{Acknowledgements}
The author is thankful to Michael Heusener, Joan Porti, Daniel Lines, Louis Funar and Rinat Kashaev for helpful remarks and suggestions related to this paper.

\bibliographystyle{amsalpha}
\bibliography{ReferencesArticleS}

\end{document}